\documentclass[reqno]{amsart}

\pdfoutput=1

\usepackage{amsfonts,amsmath,amsthm,amsxtra,amssymb,mathrsfs,dsfont}

\usepackage{tikz-cd}

\usepackage[square,numbers]{natbib}

\usepackage{graphicx}  

\usepackage{subfigure}

\usepackage{fancybox}

\usepackage[colorlinks = true, citecolor = blue, urlcolor = blue]{hyperref}

\newtheorem{theorem}{Theorem}[section]

\theoremstyle{definition}
\newtheorem{definition}[theorem]{Definition}

\newtheorem*{example}{Example}

\def\gor#1{\widetilde{#1}}
\def\tech#1{\widehat{#1}}

\def\h(#1,#2){\mbox{Hom}\left(#1,#2\right)}

\def\t(#1,#2){\mbox{Tor}\left(#1,#2\right)}
\def\e(#1,#2){\mbox{Ext}\left(#1,#2\right)}

\def\SW{\mathbb{SW}}

\def\C{\mathbb{C}}
\def\F{\mathbb{F}}

\def\N{\mathbb{N}}
\def\Q{\mathbb{Q}}
\def\R{\mathbb{R}}
\def\T{\mathbb{T}}

\def\X{\mathbb{X}}

\def\Z{\mathbb{Z}}

\def\uu{\mathbf{u}}
\def\vv{\mathbf{v}}

\def\xx{\mathbf{x}}

\def\VV{\mathbf{V}}
\def\WW{\mathbf{W}}

\def\dgm{\mathsf{dgm}}
\def\bcd{\mathsf{bcd}}

\def\birth{\mathsf{birth}}
\def\death{\mathsf{death}}

\def\<{\langle}
\def\>{\rangle}

\title{Topological Time Series Analysis}

\author[Jose Perea]{Jose A. Perea}
\thanks{This work was partially supported by the NSF under grant DMS-1622301 and  DARPA under grant HR0011-16-2-003}
\address{
\shortstack[l]{
Department of Computational Mathematics, Science \& Engineering \\
Department of Mathematics, \\
Michigan State University \\
East Lansing, MI, USA.}}
\email{joperea@msu.edu}



\begin{document}

\maketitle

\begin{abstract}
Time series are ubiquitous in our data rich world.
In what follows I will describe how ideas from dynamical systems  and topological data analysis
can be combined to gain insights from time-varying data.
We will see several applications to the live sciences and engineering,
as well as some of the theoretical underpinnings.
\end{abstract}

\section{Lorenz and the butterfly}
Imagine  you have a project involving  a crucial computer simulation.
For  an initial value $\vv_0 = (x_0,y_0,z_0) \in \R^3$,
  a sequence   $\vv_0,\ldots, \vv_n \in \R^3$ is computed in such a way that
$\vv_{j+1}$  is determined from $\vv_j$ for $j=0,\ldots, n-1$.
After the simulation is complete  you realize that a rerun is needed for further analysis.
Instead
 of initializing at $\vv_0$, which might take a while, you take a shortcut:
you input a value $\vv_j $  selected from the middle of the current results,
and the simulation runs from there   while you go for coffee.
Figure \ref{fig:Lorenz1}  is displayed on the computer monitor upon your return;
the orange curve is the sequence  $x_0,\ldots, x_n$  from the initial simulation, and the blue curve is the $x$ coordinate for the rerun initialized at $\vv_j$:
\begin{figure}[!htb]
  \centering
  \includegraphics[width=\textwidth]{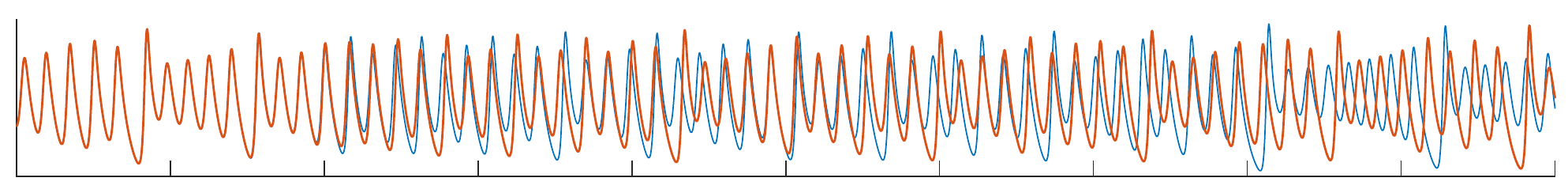}
  \caption{Orange: results from the  simulation initialized at $\vv_0$;
blue: results after manually restarting the simulation from $\vv_j$.}
  \label{fig:Lorenz1}
\end{figure}

The results agree at first, but then they diverge widely; what is going on?
 Edward Norton Lorenz,  a mathematical meteorologist, asked himself
the very same question while studying a simplified model (\ref{eq:LorenzSystem}) for weather forecasting \cite{lorenz1963deterministic}.
In the process of resolving the aforementioned discrepancy, which one could
erroneously attribute to software error or a hardware malfunction,
Lorenz laid out the foundations for what we know today as chaos theory.
The relevant set of  differential equations for the simplified model,
called the Lorenz system,  is  shown in equation (\ref{eq:LorenzSystem});
$x,y$ and $z$ are real-valued functions of time $t$,
and $\sigma, \rho, \beta \in \R$ are physical constants.
\begin{eqnarray}\label{eq:LorenzSystem}
  x'(t) &= & \sigma\cdot(y - x) \nonumber  \\
  y'(t) &= & x \cdot (\rho - z) - y\\
  z'(t)  &= & xy- \beta z \nonumber
\end{eqnarray}

Solving the Lorenz  system yields a differentiable function
\begin{equation}\label{eq:IntCurvLorenz}
\Phi: \R\times \R^3 \longrightarrow \R^3
\end{equation}
where $\Phi(t, \vv_0) = (x(t), y(t), z(t)) $ satisfies (\ref{eq:LorenzSystem}) for all $t\in \R$,  and $\Phi(0,\vv_0) = \vv_0$ for all $\vv_0 \in \R^3$.
In fact, the orange curve from Figure \ref{fig:Lorenz1} corresponds to $\Phi(t, (5,5,5))$
when  $(\sigma, \rho, \beta) = (10, 28, 8/3)$.
The discrepancy between the orange and blue curves,  as elucidated   in \cite{lorenz1963deterministic},  is a property inherent to the system.
  Lorenz realized that when manually entering $\vv_j$ as input, he only used the first few significant
digits instead of the full precision values.
In other words, the system (\ref{eq:LorenzSystem}) can be
 extremely sensitive to initial conditions
 in that  any errors are compounded exponentially with time.

\begin{figure}[!htb]
  \centering
  \includegraphics[width=0.8\textwidth]{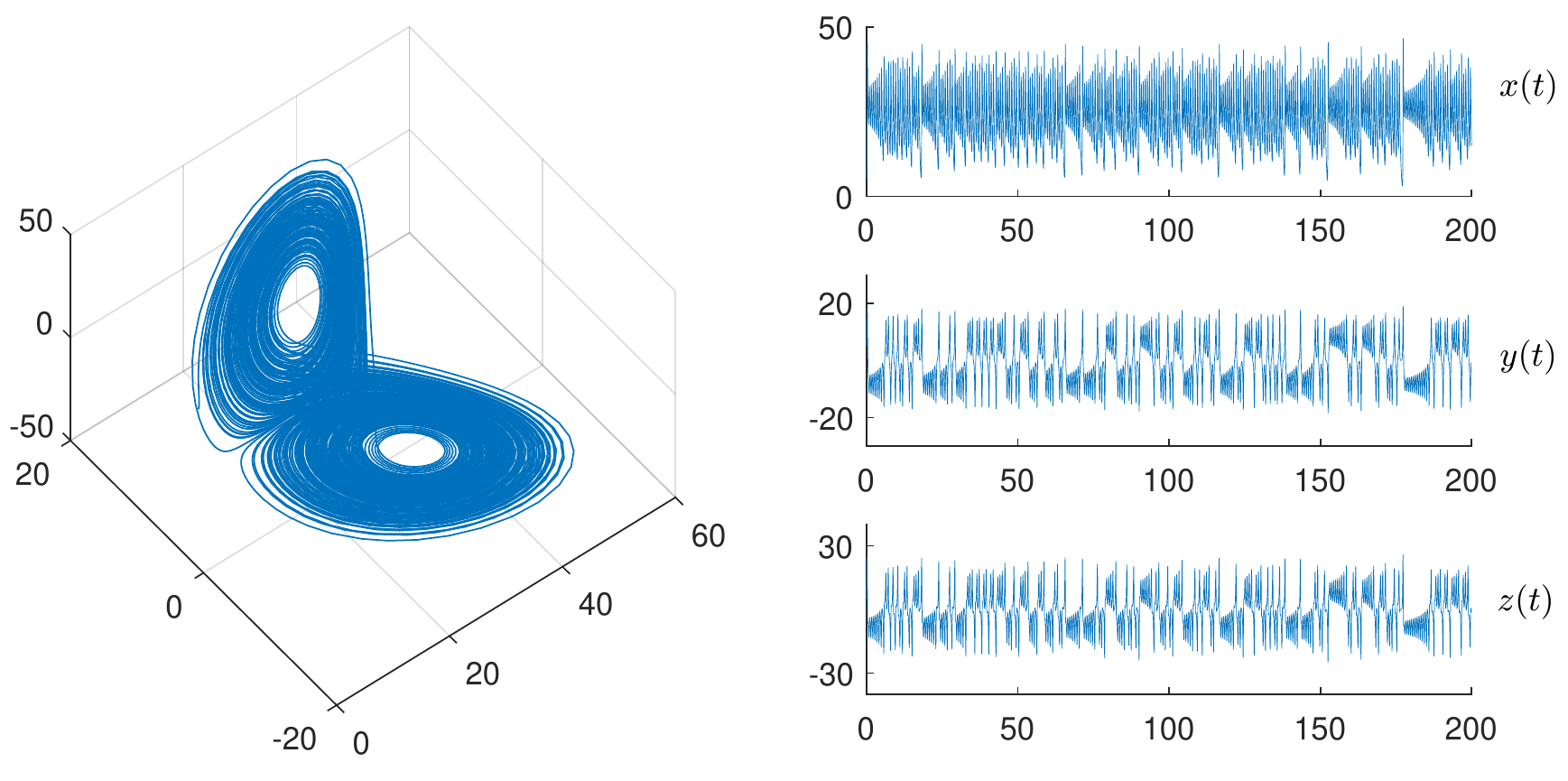}
  \caption{Lorenz's butterfly attractor}\label{fig:ButterflyAttractor}
\end{figure}

This behavior is known today as the \emph{Butterfly Effect}. The metaphor is that
 even the tiniest change in initial atmospheric conditions, e.g. such as a butterfly flapping its wings in Brazil,
may change the long-term evolution of weather patterns enough  to produce a tornado in Texas.
Such unpredictability is one of the  hallmarks of a \emph{chaotic dynamical system},
and speaks to the futility of long-term weather prediction.
The butterfly metaphor is further amplified by the shape of  the solution $t\mapsto \Phi(t, (5,5,5)) = (x(t), y(t), z(t))$,
 $0 \leq t \leq 200$,
shown in Figure \ref{fig:ButterflyAttractor} (left).

Dynamical systems are  mathematical abstractions of time-dependent physical processes.
Intuitively speaking, a  dynamical system consists of two pieces of data: a set of states $M$
--- e.g., all possible atmospheric conditions  at a given location on earth ---
along with   rules $\Phi = \{\Phi_p : p\in M\}$ describing how
each state $p\in M$ changes over time.
More specifically,
\begin{definition}
A global continuous time dynamical system is a pair $(M,\Phi)$,
where $M$ is a topological space
 and  $\Phi : \R \times M \longrightarrow M$  is a continuous map
 so that
 $\Phi(0,p) = p$, and $\Phi(s, \Phi (t, p)) = \Phi( s + t, p)$ for all $p\in M$ and all $t,s\in \R$.
\end{definition}

The typical examples arising from  differential equations  (e.g., the Lorenz system) have as state space a smooth  manifold $M$ (e.g., $\R^3$),
and the  dynamics are given by the integral curves
 (e.g., equation (\ref{eq:IntCurvLorenz}))  of a  smooth  vector field on $M$
(e.g.,  (\ref{eq:LorenzSystem})).

Some subsets of $M$ are specially important since
they attract the
  evolution of states in close proximity.
Indeed,
a set $A\subset M$   is called an attractor if it satisfies three conditions:
(1) it is compact,
(2)  it is an invariant  set  --- that is,   if  $a\in A$ then
$\Phi(t,a) \in A$ for all $t \geq 0$ ---
and
(3) it has an open basin of attraction.
In other words, there is an invariant open neighborhood $U \subset M$ of $A$,
so that
\[
\bigcap_{t\geq 0} \{\Phi(t,p) : p\in U\} = A
\]
An attractor is called  strange --- and this is just a name --- if it exhibits fractal properties.
For instance, Lorenz's butterfly $\Lambda \subset \R^3$ from Figure \ref{fig:ButterflyAttractor}
is  one of the most widely known examples of a strange attractor.

\section{Persistent Homology: Measuring  Shape  From Finite Samples}
The shape of  attractors  carries a great deal of information about the global structure of a dynamical system.
Indeed,  attractors with the shape of a circle $S^1 = \{z\in \C : |z| = 1\} $
give rise to  periodic processes;
fractal geometry  is evidence of chaotic behavior;
while  high-dimensional tori
$\mathbb{T}^n= S^1\times \cdots \times S^1$
are linked to quasiperiodicity.
The latter is a type of recurrence emerging from the superposition of periodic oscillators
with  \emph{incommensurate} (i.e., linearly independent over $\Q$) frequencies.
Quasiperiodicity  appears, for example,  in turbulent fluids    \cite{ruelle1971nature},
and the detection of biphonation  in  high-speed laryngeal  videoendoscopy \cite{tralie2018quasi} --- i.e.,
videos of vibrating vocal folds.
Similarly, the existence of chaos in brain activity is a question of considerable interest in neuroscience \cite{korn2003there},
as is the presence of periodic oscillators  in biological systems \cite{rustici2004periodic}.

A    data analysis question that has received significant attention in recent years
is how to measure the shape
of a topological space $\mathbb{X}$ --- e.g., an attractor  ---
from a finite set $X$ (the data) approximating it.
This is the type of problem driving  advances in Topological Data Analysis \cite{carlsson2014topological},
and where tools like persistent homology \cite{perea2018brief} --- which we will describe next --- are relevant.
If  $M$  is  a metric space with
  metric $\rho$, and $\X,X\subset M$, then
the quality of the approximation of  $\X$ by $X$ can be quantified
as follows.
For $\alpha \in \R$, define the $\alpha$-offset of $Z\subset M$ as
\begin{equation}\label{eq:Offset}
Z^{(\alpha)}
=
\left\{
y\in M : \inf\limits_{z\in Z}\rho(y,z) \leq \alpha
\right\}
\end{equation}
and the \emph{Hausdorff distance} between $Z$ and $\gor{Z}$ as
\begin{equation}\label{eq:HausssDist}
d_H\left(Z, \gor{Z}\right) := \inf \left\{ \alpha > 0 : Z \subset \gor{Z}^{(\alpha)} \mbox{ and } \gor{Z} \subset Z^{(\alpha)} \right\}
\end{equation}
The typical strategy in topological data analysis
is  to replace the finite sample $X$ by a space
whose shape captures that of $\X$,
provided  $d_H(X,\X)$ is small enough.
Indeed, for $Z\subset M$ and $\alpha\in \R$,
the \emph{Rips complex} of $Z$ at scale $\alpha$ is the set
\begin{equation}\label{eq:RipsFiltration}
R_\alpha(Z) := \big\{ \{z_0,\ldots, z_k\}
\subset Z:
\rho(z_i, z_j) \leq \alpha \; \mbox{ for all }\; 0 \leq i,j \leq k\big\}
\end{equation}
This is in fact a \emph{simplicial complex};
points in $Z$ can be thought of as vertices,
sets with two  elements $\{z_0, z_1\} \in R_\alpha(Z)$ are  edges,
$\{z_0,z_1,z_2\} \in R_\alpha(Z)$ is a triangular face, and so on.
A  theorem of Janko Latschev  \cite{latschev2001vietoris} contends
that if $\X$ is a closed Riemannian manifold
and
$d_H(X,\X)$ is small,  $X$ not necessarily finite,  then  (the geometric realization of)
$R_\alpha(X)$ is homotopy equivalent to
$\X$ for small  $\alpha$.

\begin{figure}[!htb]
  \centering
  \includegraphics[width=.5\textwidth]{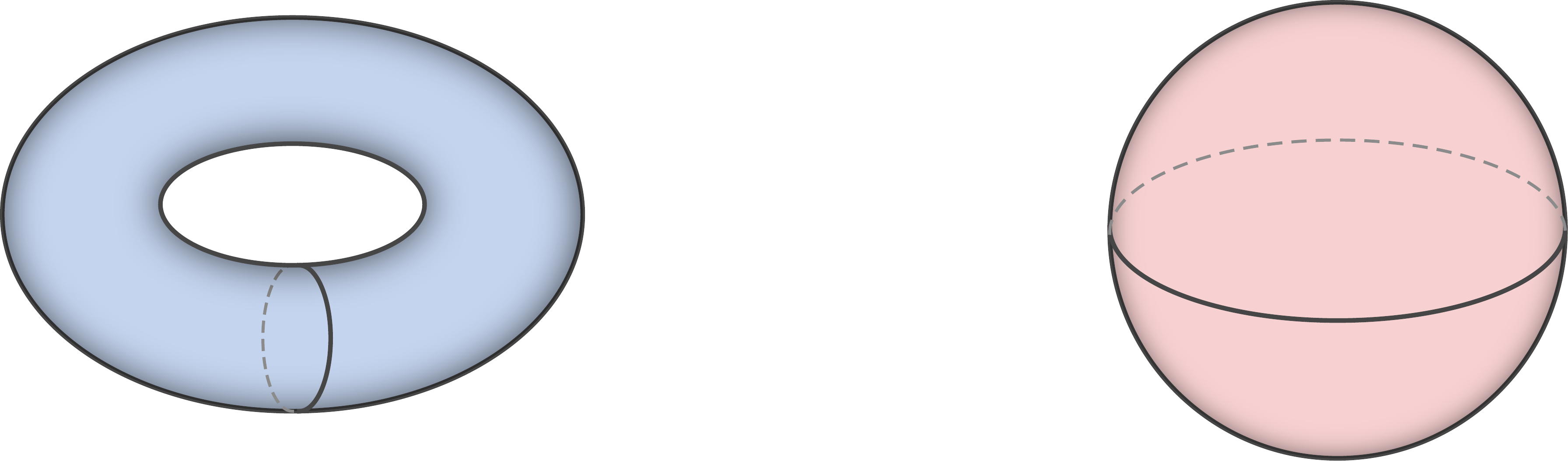}
  \caption{Left: The 2-dimensional torus $\T^2 = S^1\times S^1$.
Right: The 2-sphere $S^2 = \{\xx \in \R^3 :\|\xx\|=1\}$.}
   \label{fig:torus_circle}
\end{figure}

The homotopy type of a space refers to those properties which are invariant under continuous deformations; e.g.,
is it connected? are there holes?
and they can be quantified using  \emph{singular homology}.
Given an integer $n\geq 0$,
the $n$-th homology of a topological space $B$ with coefficients in a field $\F$ is a vector space $H_n(B;\F)$.
Its dimension  $\beta_n(B;\F)$ ---
 the
$n$-th \emph{Betti number} of $B$ with coefficients in $\F$ ---
provides a count for the number of essentially distinct $n$-dimensional holes in $B$.
Indeed, if $\F = \Q$, then $\beta_0(B) = \beta_0(B;\Q)$ counts the number of path-connected components of $B$,
$\beta_1(B)$ is the number of essentially distinct loops in $B$ bounding a hole, $\beta_2 (B)$
is roughly the number of closed 2-dimensional regions bounding a void, and so on for $n\geq 3$.
Here is an example: the 2-dimensional torus $\T^2 = S^1\times S^1$ in Figure \ref{fig:torus_circle} (left)
has Betti numbers $\beta_0(\T^2) = 1$ since it is path-connected,
$\beta_1(\T^2) =2 $ since it has a   horizontal and a vertical hole,
$\beta_2(\T^2) =1$ since $\T^2$ encloses an empty volume, and $\beta_n(\T^2) = 0$ for all $n \geq 3$.
Similarly, the 2-sphere $S^2  = \{\xx\in \R^3 : \|\xx\|=1\}$   has Betti numbers
$\beta_0(S^2) = \beta_2(S^2) = 1$, $\beta_1(S^2)= 0$,  and $\beta_n(S^2) =0 $ for  $n\geq 3$.

Continuous maps $f: B \longrightarrow B' $ between topological spaces
can also be studied through the lens of homology.
Indeed, any such $f$
induces   linear transformations $f_n : H_n(B;\F) \longrightarrow H_n(B';\F)$, $n\geq 0$,
between the corresponding homology vector spaces
 in such a way that if $id : B \longrightarrow B$ is the identity
of $B$, then $id_n$ is the identity of $H_n(B;\F)$, and $(f\circ g)_n = f_n \circ g_n$ whenever compositions make sense.

%

For realistic data $X\subset M$,   the Betti numbers   $ \alpha \mapsto \beta_n\left(R_{\alpha}(X); \F\right)$
are expected to be fairly unstable as  $\alpha$ varies.
Indeed, sampling artifacts or noise in $X$ can produce holes which are present in $R_{\alpha}(X)$
but not in $R_{\alpha + \delta}(X)$ for small $\delta > 0$
(e.g., see Figure \ref{fig:BarcodeRips_Circle} below).
This is where  comparing the homology of spaces related by maps is useful.
Indeed, let $\mathcal{Z} = \{Z_\alpha\}_{\alpha \in \R}$ be a family of spaces
such that $Z_\alpha \subset Z_{\alpha'}$ for all pairs $\alpha \leq \alpha'$.
Such  a family is called a \emph{filtration}, and
we have already seen two examples:
The offset filtration $\left\{Z^{(\alpha)}\right\}_{\alpha \in \R}$  defined in equation (\ref{eq:Offset}),
and   the Rips filtration $\mathcal{R}(Z) = \{ R_\alpha(Z)\}_{\alpha \in \R}$
from (\ref{eq:RipsFiltration}).
Given a filtration  $\mathcal{Z}= \{Z_\alpha\}_{\alpha \in \R}$ and $\alpha \leq \alpha'$,
the inclusion map $\iota^{\alpha,\alpha'} : Z_\alpha \hookrightarrow Z_{\alpha'}$ induces linear transformations
\begin{equation}\label{eq:PersistentHomology}
\iota^{\alpha,\alpha'}_n : H_n(Z_{\alpha};\F) \longrightarrow H_n(Z_{\alpha'};\F)
\;\;\; , \;\; n\geq 0
\end{equation}
Classes in $H_n(Z_\alpha;\F)$ which are not in the Kernel of $\iota_n^{\alpha,\alpha'}$
 for large $\alpha' - \alpha$, are interpreted as being persistent throughout the filtration.
The  collection of vector spaces and linear maps resulting form (\ref{eq:PersistentHomology}) is called
the $n$-dimensional \emph{persistent homology}, with coefficients in $\F$, of the filtration $\mathcal{Z}$.
More generally,
\begin{definition}
A persistence vector space $\VV$ is a collection of  vector spaces $ V_{\alpha}$, $\alpha \in \R$,
and  linear transformations $\iota^{\alpha, \alpha'} : V_\alpha \longrightarrow V_{\alpha'}$,
  $\alpha \leq \alpha'$, so that:
\begin{enumerate}
\item $\iota^{\alpha,\alpha}$ is the identity of $V_\alpha$ for every $\alpha \in \R$.
\item $\iota^{\alpha', \alpha''} \circ \iota^{\alpha, \alpha'} = \iota^{\alpha, \alpha''}$, whenever $\alpha \leq \alpha' \leq \alpha''$.
\end{enumerate}
Two persistence vector spaces $\VV =  \{V_\alpha, \iota^{\alpha, \alpha'} \}$
and
$\WW = \{W_\alpha, \kappa^{\alpha, \alpha'}\}$
are   isomorphic, denoted $\VV \cong \WW$,  if there are linear isomorphisms
$ T_\alpha : V_\alpha \longrightarrow W_\alpha$ for all $\alpha \in \R$,
so that $\kappa^{\alpha,\alpha'} \circ T_\alpha = T_{\alpha'} \circ \iota^{\alpha,\alpha'}$
whenever $\alpha \leq \alpha'$.
\end{definition}

We will concentrate on three quantities  for a nonzero element $\gamma \in V_\alpha$:
\begin{eqnarray}
\birth(\gamma) &:=&
\inf
\left\{
\gor{\alpha} \leq \alpha :
\gamma \in \mathsf{Image}\left(\iota^{\gor{\alpha},\alpha}\right)
\right\}
\nonumber\\[.1cm]
\death(\gamma) &:=&
\sup
\left\{
\alpha' \geq \alpha :
\gamma \notin \mathsf{Kernel}\left(\iota^{\alpha, \alpha'}\right)
\right\}
\\[.1cm]
\mathsf{persistence}(\gamma) &:= &  \death(\gamma) - \birth(\gamma) \nonumber
\end{eqnarray}
As we saw before, the Betti numbers capture the homology of a space and hence yield a succinct
shape descriptor for its  topology.
Persistent homology, on the other hand, describes the evolution of homological features in a filtration.
Under favorable circumstances    \cite{crawley2015decomposition} it can  be   described
by an invariant  called \emph{the barcode}:
\begin{theorem}
Let $\VV =  \{V_\alpha, \iota^{\alpha,\alpha'}\}$ be a persistence vector space so that
$\dim(V_\alpha)$ is finite for all $\alpha \in \R$.
Then, there  exists  a  multiset
(i.e., a set whose elements may appear with repetitions) of intervals
$ I \subset [-\infty, \infty] $
 called the barcode of $\VV$, and denoted $\bcd(\VV)$,   so that:
\begin{enumerate}
\item It subsumes the Betti numbers: If $\alpha\in \R$, then $\dim(V_\alpha)$ is exactly the number of intervals
$I\in \bcd(\VV)$, counted with repetitions, so that $\alpha\in I$.
\item It encodes persistence: For every $I \in \bcd(\VV)$ and every $\alpha \in I$, there exists
  $\gamma \in V_\alpha$ so that the left and right end-points of $I$ are
$\birth(\gamma)$ and  $\death(\gamma)$, respectively.
\item It is an invariant: If $\WW$ is a  persistence vector space with $\dim(W_\alpha) < \infty$ for all $\alpha \in \R$,
then $\bcd(\VV) = \bcd(\WW)$   if and only if
$\VV \cong \WW$.
\end{enumerate}
\end{theorem}

We will use $\bcd_n^\mathcal{R}(X;\F)$ to denote
the barcode for the  $n$-dimensional persistent homology of the Rips filtration on a metric space $X$.
Below in Figure \ref{fig:BarcodeRips_Circle} we show an example  for
$X\subset \R^2$ sampled with noise around the unit circle
$S^1$,  the rips complex $R_\alpha(X)$ at $\alpha = 0$,   $ 0.36 $, $0.6$, $1.21$,
and  the intervals (i.e., the horizontal blue lines) which comprise the barcode $\bcd^\mathcal{R}_1(X;\Z/2)$.
The computations were performed using the C++ library  \verb"Ripser" \cite{bauer2017ripser}.
The long interval is indicative of a persistent 1-dimensional hole in the data,
which is consistent with $X$ being sampled
around $S^1$; recall that $\beta_1(S^1;\Z/2)= 1$.
The short-lived intervals, on the other hand, are due to noise and sampling artifacts.

\begin{figure}[!htb]
  \centering
  \includegraphics[width=0.9\textwidth]{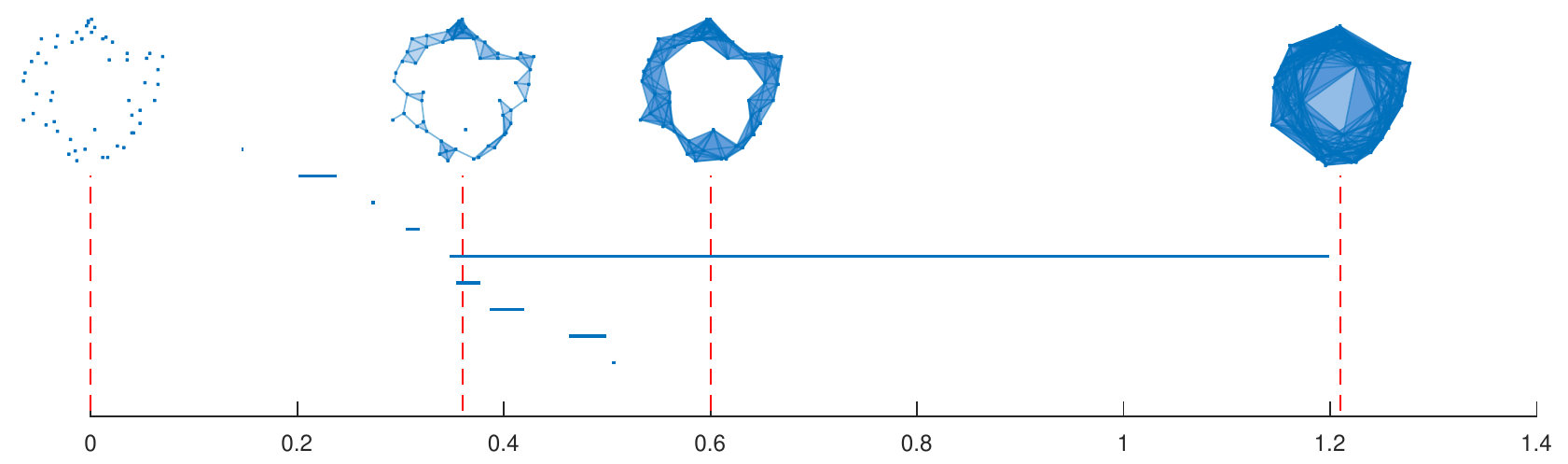}
  \caption{Barcode for the Rips filtration on $X\subset \R^2$ near $S^1$.}\label{fig:BarcodeRips_Circle}
\end{figure}

As we have seen, persistent homology can be used to infer the topology
of  a space  $\X$ given a finite sample $X$:
the number of long intervals in $\bcd_n^\mathcal{R}(X;\F)$ suggests
a value for $\beta_n(\X;\F)$.
The barcodes $\bcd_n^\mathcal{R}(X;\F)$ can also be used to quantify/indetify
properties of dynamical systems: if $\X$ is an attractor,
a barcode like the one in Figure \ref{fig:BarcodeRips_Circle} would
be indicative of periodicity, while the barcodes
in Figure \ref{fig:SWP_example} would point to  quasiperiodicity.
Moreover, measures of  fractal geometry can also be derived  from persistent homology;
see for instance \cite{macpherson2012measuring} and \cite{robins1999towards}.

\section{Attractor Reconstruction: Time series data,  Takens'  Theorem and Sliding Window Embeddings}
In practice it is exceedingly  rare to have an explicit  mathematical  description of a dynamical system
 of interest.
Instead,  often times one can gather measurements of relevant quantities  for each sate $p\in M$
--- e.g.,
in weather prediction one can estimate  temperature,  pressure, etc.
A  way of measuring can be thought of as a  continuous map
$F : M \longrightarrow \R$ --- called an \emph{observation function} ---
and   given an initial state $p \in M$,
one obtains the time series
\begin{equation}\label{eq:TseriesFromDynamics}
\begin{array}{rccl}
\varphi_p :& \R &  \longrightarrow & \R \\
  & t & \mapsto &  F \circ\Phi(t,p)
\end{array}
\end{equation}
For example, the blue and orange curves from Figure \ref{fig:Lorenz1}
are examples of time series from the Lorenz system.
A single time series may appear to be a complete oversimplification
of the underlying dynamics.
The Takens' embedding theorem, due to Floris Takens \cite{takens1981detecting},
implies  that they can actually be very useful:

\begin{theorem}
Let $M$ be a smooth, compact, Riemannian manifold; let $\tau > 0$ be a real number  and let $d \geq 2\dim(M)$ be an integer.
If  $\Phi \in C^2(\R \times M,\, M)$ and $F \in C^2(M, \R)$ are generic,
then for $\varphi_p(t)$ defined by (\ref{eq:TseriesFromDynamics}),
the delay map
\begin{equation}\label{eq:TakenEmbd}
\begin{array}{rccl}
  \varphi : & M & \longrightarrow  &  \R^{d + 1} \\
   & p & \mapsto &  \big(\varphi_p(0),  \varphi_p(\tau), \varphi_p(2\tau) , \ldots ,  \varphi_p(d\tau) \big)
\end{array}
\end{equation}
is an embedding.
\end{theorem}
Here generic means that the set of functions $ \Phi, F $ for which
(\ref{eq:TakenEmbd}) is an embedding
forms an open and dense set with respect to the Whitney topology.
In fact, if $A\subset M$ is a strange attractor, then
$\varphi$ restricted to $ A$ will be (generically) an embedding whenever
$d  $ is at least twice  the box counting (a notion of fractal) dimension of $A$.
 Takens'  theorem motivates the following definition,

\begin{definition}
Let $ f : \R \longrightarrow \R$ be a function, $\tau > 0$ a real number and $d > 0$  an integer.
The \emph{sliding window embedding} of $f$, with parameters $d$ and $\tau$, is the vector-valued function
\begin{equation}
\begin{array}{rccl}
SW_{d,\tau} f : &  \R & \longrightarrow & \R^{d+1 } \\
&t& \mapsto & \big(f(t), f(t + \tau), f(t + 2\tau),  \ldots , f(t + d\tau) \big)
\end{array}
\end{equation}
The integer $d+1$ is the dimension, $\tau$ is the delay and the product $d\tau$ is the window size.
For $T\subset \R$, the set
\begin{equation}\label{eq:SWPointCloud}
\SW_{d,\tau}f = \{SW_{d,\tau} f(t) : t\in T\}
\end{equation} is the \emph{sliding window point cloud} associated
to $T$.
\end{definition}

\begin{figure}[!htb]
  \centering
  \includegraphics[width= 0.9\textwidth]{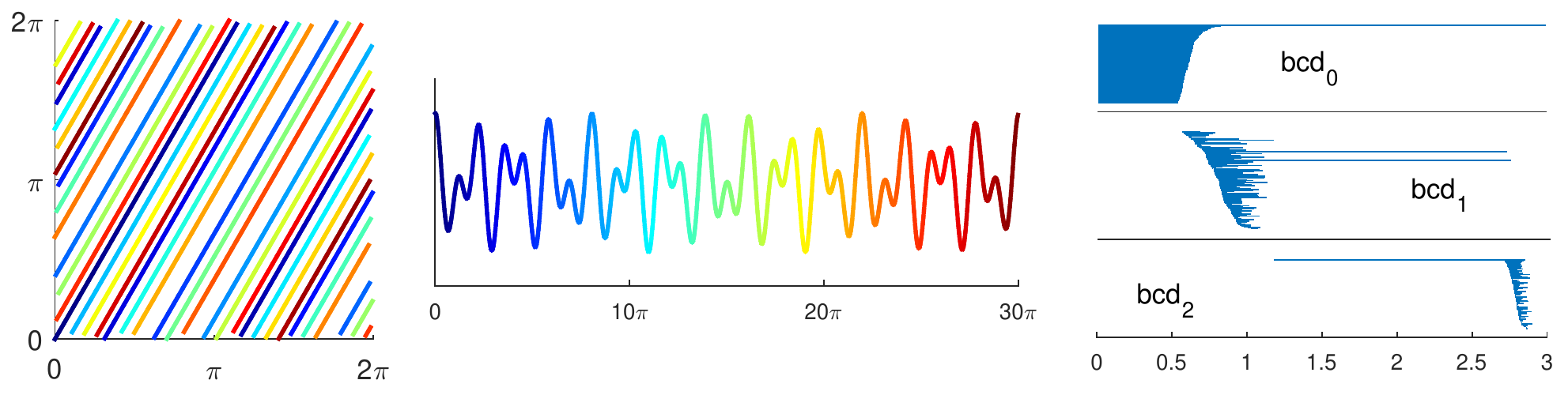}
  \caption
{Left: The dynamics $\Phi$ on the torus.
 The colors, blue through red, indicate the time $t$
for  $t\mapsto \Phi\big(t, (1,1)\big)\in \T^2$.
Center: The time series  $f(t) = F\circ \Phi(t, (1,1)) = \cos(t) + \cos(\sqrt{3}t)$.
The colors indicate the time variable $t \in [0, 30\pi]$ and are coordinated
 with the torus on the left panel. Right:
Barcodes for the Rips filtration  $\mathcal{R}\big(\SW_{d,\tau} f\big)$;
the number of long intervals recovers the Betti numbers of
the attractor:
$\beta_0(\T^2) = \beta_2(\T^2)=1, \beta_1(\T^2) =2 $.}
\label{fig:SWP_example}
\end{figure}

Hence,
given time series data
$f(t) = \varphi_p(t)$
observed from a potentially unknown  dynamical system
$(M,\Phi)$, Takens' theorem implies that (generically) the
sliding window point cloud
$\SW_{d,\tau} f$ provides a topological copy   of
$\{\Phi(t,p) : t\in T\}\subset M$.
In particular, this will reconstruct attractors.
The underlying shape of $\SW_{d,\tau}f$ can then be quantified with persistent homology,
and the associated barcodes  $\bcd_n^\mathcal{R}(\SW_{d,\tau}f; \F)$
can be used as features in inference, classification,
and learning tasks  \cite{garland2016exploring, robins1999towards, xu2018twisty}.
We will see shortly  several applications of these ideas to science and engineering;
for now, here is an instantiation of the pipeline:

\begin{example}
Let $\omega\in \R$ be  irrational;
we will use $\omega = \sqrt{3}$ for computations
but any other choice would do.
Consider the  dynamics $\Phi $ and the observation function $F$
on the torus $\T^2 = S^1\times S^1 \subset \C^2$, given by
\[
\begin{array}{rcclcrccl}
\Phi : &\R \times \T^2& \longrightarrow& \T^2
&&
F : & \T^2 & \longrightarrow & \R
\\
&\big(t, (z_1, z_2)\big)& \mapsto& \left(e^{ it}z_1,  e^{  i\omega t}z_2\right)
&&&
(z_1,z_2) &\mapsto & \mathsf{Re}(z_1 + z_2)
\end{array}
\]
If $p   \in \T^2$, then   $\{\Phi(t,p) : t\in \R\}$ is dense in $\T^2$,
  and hence $\T^2$ is the only attractor;
e.g., see  \cite{boothby1986introduction}, page 86,
Example 6.15.
For $p = (1,1)$ we obtain the quasiperiodic time series
$ f (t) = F \circ \Phi\big(t,(1,1)\big)= \cos( t) + \cos(\omega t)$,
and we show  in Figure \ref{fig:SWP_example}
the dynamics $t\mapsto \Phi(t,(1,1))$ on the torus (left),
the resulting time series $f(t)$ (center), and the
barcodes $\bcd_n = \bcd_n^\mathcal{R}(\SW_{d,\tau} f; \Z/2)$
in dimensions $n = 0,1,2$ (right).
\end{example}

\section{Some Applications of Sliding Window Persistence}

\subsection*{Wheeze Detection} A wheeze is an abnormal whistling  sound produced while breathing.
It is often associated with obstructed airways and lung diseases such as
 asthma, lung cancer and congestive heart failure.
In \cite{emrani2014persistent},
 Emrani,   Gentimis and   Krim show  that the 1-dimensional
 barcode $\bcd_1^\mathcal{R}(\SW_{d,\tau} f ; \F)$,
and particularly the length of its longest interval
(i.e., its maximum persistence)
\begin{equation}\label{eq:MaxPers}
    \mathsf{mp}_1^\mathcal{R}(\SW_{d,\tau} f ; \F)
       =
    \max\left\{\mathsf{length}(I) : I \in \bcd_1^\mathcal{R}(\SW_{d,\tau} f ;\F )  \right\}
\end{equation}
is an effective feature for wheezing detection in  recorded breathing
sounds.
Indeed, the presence of wheezing   leads to  circular sliding window point clouds.
When testing on a large database of sound recordings, Emrani et. al. shown that
(\ref{eq:MaxPers}) leads to a higher detection accuracy  than that     of competing methods.

\subsection*{Periodicity Quantification in Gene Expression Data}
Many biological processes, including the cell cycle, cell division,  and the circadian clock,
are periodic in nature.
An important problem in systems biology is to describe these mechanisms  at a  genetic level \cite{rustici2004periodic}.
Biologists approach this   by first collecting data;
specifically, how gene expression
changes across time in a given model organism --- e.g., yeast, mice, etc.
Given said data,
it is shown  in \cite{perea2015sw1pers}
that  measures similar to (\ref{eq:MaxPers}) can
outperform  state-of-the-art methods  for periodicity quantification,
and lead to the discovery of novel clock-regulated genes.

\subsection*{Segmentation of Dynamic Regimes}
Complex  high-dimensional  systems can  exhibit abrupt changes in qualitative behavior.
For instance, the hearth's climate has undergone
several  sudden transitions to and back from  a ``snowball earth''  \cite{hoffman1998neoproterozoic}.
Identifying markedly different regimes in  a system's evolution
can thus be used for warning, modeling and parameter estimation purposes
\cite{dakos2008slowing}.
Berwald et. al. show in \cite{berwakld20014automatic} that, given time series data,
  effective classifiers can be trained on features from the barcodes of sliding window point clouds,
with the goal of automatically segmenting a system into different behavioral regimes.
Some applications of their methodology include the detection of bifurcations
on stochastic and chaotic systems, as well as the analysis of temperature and CO$_2$ levels in ice cores.

\subsection*{Chatter Detection and Classification in Machining}
Turning and milling are cutting processes used extensively in industrial manufacturing.
Chatter, or machining vibrations, are wide oscillations of the cutting
tool with respect to the metal workpiece; these undesired movements leave surface flaws on the production piece during turning and milling.
Khasawneh et. al. show in  \cite{khasawneh2016chatter, khasawneh2018chatter}
that, given time series describing the undulations of the cutting piece,
$\bcd^\mathcal{R}_1(\SW_{d,\tau} f; \F)$ can be used as input features in classification algorithms for chatter detection.

\subsection*{Motion Recognition from Wearable Sensors}
Dirafzoon et. al. show in \cite{dirafzoon2016action} that
 barcodes on sliding window point clouds from motion capture time series,
can be used to classify activities with very high accuracy.

\subsection*{Periodicity and Quasiperiodicity in Video Data}
A video can be thought of as sampling a function $f: \R \longrightarrow \R^k$.
It follows that computing $\bcd^\mathcal{R}_n(\SW_{d,\tau} f ;\F)$ for the sliding window point
$\SW_{d,\tau} f \subset \R^{k(d+1)}$
can be used to detect recurrent behavior in video data, without the need for tracking or
surrogate signals.
See for instance
\cite{tralie2018quasi}
for  applications of these ideas to the problem of quantifying periodicity
and quasiperiodicity in video data,
as well as the synthesis of
slow-motion  videos from recurrent movements \cite{tralie2018topological}.

%
%

\section{Theoretical Investigations of Sliding Window Persistence}
The barcodes $\bcd_n^\mathcal{R}(\SW_{d,\tau} f; \F)$ have been used  successfully  in several  applications.
However, very little is known about  their theoretical  behavior with respect to
the parameters $d,\tau $, the field $ \F$ and the homological dimension $n$.
One of the main difficulties lies in that our knowledge of
the homotopy type/homology of the Rips complex $R_\alpha(X)$, for arbitrary $\alpha$,
is rather limited:
planar circles are essentially the only spaces where we have complete answers \cite{adamaszek2017vietorisCircle}.
I would like to highlight  next some of what we do know.

\subsection*{Sliding Window Persistence of Periodic Functions}
Let $\T = \R/2\pi\Z$   and let
$  C^1(\T, \R)$
denote the set of continuously differentiable real-valued functions
on $\T $.
As a warm up example, let   $\zeta(t) = \sin(Lt + \phi)$, for $L\in \N$ and $\phi \in \R$.
A bit of trigonometry shows that
\[
SW_{d,\tau} \zeta(t) = \sin(Lt + \phi)\uu + \cos(Lt  + \phi)\vv
\]
where $\uu = SW_{d,\tau}\cos(Lt )|_{t = 0}$ and $\vv = SW_{d,\tau} \sin(L t )|_{t = 0}$ .
It readily follows that
if $d\geq 1$ and  the set $\{\uu,\vv\}$ is linearly independent,
then
$\SW_{d,\tau} \zeta = SW_{d,\tau} \zeta(\T)$ is
 a planar ellipse.
The semi-major and semi-minor axes
can be computed explicitly as
\begin{equation}\label{eq:SemiAxes}
a =
\sqrt{
\frac
{(d + 1) + \left|\frac{\sin(L(d+1)\tau)}{\sin(L\tau)}\right|}
{2}
}
\;\;\;\;\; \mbox{ and }\;\;\;\;\;
b = \sqrt{\frac{(d + 1) - \left|\frac{\sin(L(d+1)\tau)}{\sin(L\tau)}\right|}{2}}
\end{equation}

The persistent homology of the Rips filtration on
ellipses with small eccentricity, i.e. when $b < a < \sqrt{2}b$,
has been recently studied by
Adamaszek et. al. \cite{adamaszek2017vietoris}.
In particular, their work implies that  if
\[
\alpha_1 = \frac{ 4\sqrt{3}ab^2}{ a^2 + 3b^2}
\;\;\;\;\; \mbox{ and }\;\;\;\;\;
\alpha_2 = \frac{4\sqrt{3}a^2b}{3a^2 + b^2}
\]
then the homotopy type of $R_\alpha(\SW_{d,\tau} \zeta)$
for $ 0 < \alpha \leq \alpha_2$ is either that of the circle, or that of a wedge
of 2-dimensional spheres as follows:
\[
R_\alpha(\SW_{d,\tau} \zeta)
\simeq
\left\{
  \begin{array}{lcl}
    S^1 & \hbox{ for } & 0 < \alpha < \alpha_1 \\
    S^2 & \hbox{ for } &\alpha = \alpha_1 \\
    \bigvee^5 S^2 & \hbox{ for }& \alpha_1 < \alpha  < \alpha_2 \\
    \bigvee^3 S^2 & \hbox{ for } &\alpha = \alpha_2
   \end{array}
\right.
\]
The range  $\alpha_1 < \alpha < \alpha_2$ is specially  interesting since
 only one of the five
2-dimensional classes   persists.
In other words,  the linear transformation
\[
H_2(R_\alpha(\SW_{d,\tau} \zeta) ;\F)
\longrightarrow
H_2(R_{\alpha'}(\SW_{d,\tau} \zeta) ;\F)
\]
has rank one for every $\alpha_1 < \alpha < \alpha' < \alpha_2$.
We readily obtain the following,

\begin{theorem} Let $\zeta(t) = \sin(Lt + \phi)$
for $L\in \N$, $\phi \in \R$, and for $d\in \N$, $\tau > 0$
let $a,b \geq 0$
be as in (\ref{eq:SemiAxes}).
If
$a < \sqrt{2}b$ and $\F$ is a field,
then the maximum persistence (\ref{eq:MaxPers}) in dimensions 1 and 2 satisfies
\begin{eqnarray*}
\mathsf{mp}^\mathcal{R}_1( \SW_{d,\tau} \zeta; \F) &=& \frac{ 4\sqrt{3}ab^2}{ a^2 + 3b^2} \\
\mathsf{mp}^\mathcal{R}_2( \SW_{d,\tau} \zeta; \F)
&\geq&
\frac{4\sqrt{3}a^2b}{3a^2 + b^2} - \frac{ 4\sqrt{3}ab^2}{ a^2 + 3b^2} \\
\end{eqnarray*}
\end{theorem}

The cases $\alpha > \alpha_2$  for $ a < \sqrt{2}b$,
and $\alpha > 0$ for $a \geq \sqrt{2}b$, are currently open.
The case $a = b$, i.e.
when $\SW_{d,\tau} \zeta$ is a circle,  is much better understood.
This happens  when
the window size $d\tau$ is equal to (an integer multiple of) $\frac{d}{d+1} \frac{2\pi}{L}$,
which is a little bit under   $\mathsf{Period}(\zeta) = \frac{2\pi}{L}$.
Since the homotopy type of $R_\alpha(S^1)$
is known for all $\alpha > 0$ \cite{adamaszek2017vietorisCircle},  we get that
\[
R_\alpha(\SW_{d,\tau} \zeta)
\simeq
\left\{
  \begin{array}{lll}
    S^{2k + 1} & \hbox{ for } & \sin\left(\frac{\pi k}{2k + 1}\right) < \frac{\alpha}{\sqrt{2(d+1)}} <  \sin\left(\frac{\pi (k+1)}{2k + 3}\right)\\ \\
    \bigvee^{|\R|}S^{2k} & \hbox{ for } & \alpha = \sqrt{2(d+1)}\sin\left(\frac{\pi k}{2k + 1}\right)
  \end{array}
\right.
\]
where the linear transformation
\[
H_{2k+1}(R_\alpha(\SW_{d,\tau}\zeta);\F)
\longrightarrow
H_{2k+1}(R_{\alpha'}(\SW_{d,\tau}\zeta);\F)
\]
is an isomorphism (with rank one) for every
\[
\sqrt{2(d+1)}\sin\left(\frac{\pi k}{2k + 1}\right)
<
\alpha \leq \alpha'
<
\sqrt{2(d+1)}\sin\left(\frac{\pi (k+1)}{2k + 3}\right).
\]
Therefore,
\begin{theorem} Let $\zeta(t) = \sin(Lt + \phi)$,
let $\F$ be a field,
and let $\tau = \frac{2\pi}{L(d+1)}$.
Then, for every   integer  $k \geq 1$,
$\mathsf{mp}_{2k}^\mathcal{R}
(\SW_{d,\tau} \zeta; \F) = 0 $ and
\[
\mathsf{mp}_{2k-1}^\mathcal{R}
(\SW_{d,\tau} \zeta; \F)
=
\sqrt{2(d+1)}\cdot
\left(
\sin\left(\frac{\pi k}{2k + 1}\right)
-
\sin\left(\frac{\pi (k-1)}{2k - 1}\right)
\right)
\]
\end{theorem}


One strategy to understand
$\dgm_n^\mathcal{R}(\SW_{d,\tau} f ; \F)$ for a
general function $f\in C^1(\T,\R)$, is to first
approximate $f$ by its truncated Fourier series
\[
S_N f (t) = \sum_{|n| \leq N} \tech{f}(n) e^{ i n t}
\;\;\; , \;\;\;
\tech{f}(n) = \frac{1}{2\pi}\int_\T f(t) e^{-i nt} d t
\]
and then investigate the asymptotic behavior of
the sequence  of barcodes
\begin{equation}\label{eq:SW_TruncatedFourier}
\bcd_n^\mathcal{R}(\SW_{d,\tau} S_N f ; \F)
\;\;\;\; , \;\;\; N \in \N
\end{equation}
Indeed, the  analysis of
$\zeta(t) = \sin(Lt + \phi)$
presented earlier
can be bootstrapped to
 trigonometric polynomials,
and the stability Theorem \cite{chazal2009proximity}
can be used to
study $\bcd_n^\mathcal{R}(\SW_{d,\tau} f; \F)$
via   the  behavior of
(\ref{eq:SW_TruncatedFourier}) as $N \to \infty$.
This line of reasoning was explored in \cite{perea2015sliding}.
In particular, it results in sights for the choice of window size
(it should be approximately  $\frac{d}{d+1}$ times the period length),
the embedding dimension (larger than twice the number of
relevant harmonics) and the choice of
field of coefficients (one whose characteristic does not
divide $\frac{2\pi}{\mathsf{Period}(f)}$).
We end with a theorem \cite[6.8]{perea2015sliding} relating the sliding window
persistence of $f$ to its harmonic content:

\begin{theorem}
 Let $f\in C^1(\T, \R)$ be so that
$f(t + \frac{2\pi}{L}) = f(t)$ for all $t\in \T$,
and assume (for simplicity) that $f$ has been centered and normalized:
\[
\tech{f}(0) = 0
\;\;\;\; , \;\;\;\;
\|f\|_2 : =  \left(\frac{1}{2\pi}\int_\T |f(t)|^2 dt \right)^{1/2 }= 1
\]
Let $\tau_d =  \frac{2\pi}{L(d+1)}$,
$T\subset \T$ and
$\SW_{d,\tau_d} f = SW_{d,\tau_d}f(T)$.
If $f'(t) = \frac{df}{dt}$, then
\[
\sup_{d \in \N} \mathsf{mp}_1^\mathcal{R}(\SW_{d,\tau_d} f ; \Q)
\;\geq\;
2\sqrt{3}\sup_{n\in \Z}
\left|\tech{f}(n)\right|
- 2\sqrt{2}\left\|f'\right\|_2
d_H(T,\T)
\]
\end{theorem}

\subsection*{Beyond Periodicity}
Very little is known about the sliding window persistence of other families of functions;
there are a few results for quasiperiodic functions
\cite{perea2016persistent},
but the rest of the landscape  is essentially uncharted territory.


\begin{thebibliography}{30}
\providecommand{\natexlab}[1]{#1}
\providecommand{\url}[1]{\texttt{#1}}
\expandafter\ifx\csname urlstyle\endcsname\relax
  \providecommand{\doi}[1]{doi: #1}\else
  \providecommand{\doi}{doi: \begingroup \urlstyle{rm}\Url}\fi

\bibitem[Adamaszek and Adams(2017)]{adamaszek2017vietorisCircle}
M.~Adamaszek and H.~Adams.
\newblock The vietoris--rips complexes of a circle.
\newblock \emph{Pacific Journal of Mathematics}, 290\penalty0 (1):\penalty0
  1--40, 2017.

\bibitem[Adamaszek et~al.()Adamaszek, Adams, and Reddy]{adamaszek2017vietoris}
M.~Adamaszek, H.~Adams, and S.~Reddy.
\newblock On vietoris--rips complexes of ellipses.
\newblock \emph{Journal of Topology and Analysis}, pages 1--30.

\bibitem[Bauer(2017)]{bauer2017ripser}
U.~Bauer.
\newblock Ripser: a lean c++ code for the computation of vietoris-rips
  persistence barcodes.
\newblock 2017.
\newblock Software available at \url{https://github.com/Ripser/ripser}.

\bibitem[Berwald et~al.(2014)Berwald, Gidea, and
  Vejdemo-Johansson]{berwakld20014automatic}
J.~J. Berwald, M.~Gidea, and M.~Vejdemo-Johansson.
\newblock Automatic recognition and tagging of topologically different regimes
  in dynamical systems.
\newblock \emph{Discontinuity, Nonlinearity, and Complexity}, 3\penalty0
  (4):\penalty0 413--126, 2014.

\bibitem[Boothby(1986)]{boothby1986introduction}
W.~M. Boothby.
\newblock \emph{An introduction to differentiable manifolds and Riemannian
  geometry}, volume 120.
\newblock Academic press, 1986.

\bibitem[Carlsson(2014)]{carlsson2014topological}
G.~Carlsson.
\newblock Topological pattern recognition for point cloud data.
\newblock \emph{Acta Numerica}, 23:\penalty0 289, 2014.

\bibitem[Chazal et~al.(2009)Chazal, Cohen-Steiner, Glisse, Guibas, and
  Oudot]{chazal2009proximity}
F.~Chazal, D.~Cohen-Steiner, M.~Glisse, L.~J. Guibas, and S.~Y. Oudot.
\newblock Proximity of persistence modules and their diagrams.
\newblock In \emph{Proceedings of the twenty-fifth annual symposium on
  Computational geometry}, pages 237--246. ACM, 2009.

\bibitem[Crawley-Boevey(2015)]{crawley2015decomposition}
W.~Crawley-Boevey.
\newblock Decomposition of pointwise finite-dimensional persistence modules.
\newblock \emph{Journal of Algebra and its Applications}, 14\penalty0
  (05):\penalty0 1550066, 2015.

\bibitem[Dakos et~al.(2008)Dakos, Scheffer, van Nes, Brovkin, Petoukhov, and
  Held]{dakos2008slowing}
V.~Dakos, M.~Scheffer, E.~H. van Nes, V.~Brovkin, V.~Petoukhov, and H.~Held.
\newblock Slowing down as an early warning signal for abrupt climate change.
\newblock \emph{Proceedings of the National Academy of Sciences}, 105\penalty0
  (38):\penalty0 14308--14312, 2008.

\bibitem[Dirafzoon et~al.(2016)Dirafzoon, Lokare, and
  Lobaton]{dirafzoon2016action}
A.~Dirafzoon, N.~Lokare, and E.~Lobaton.
\newblock Action classification from motion capture data using topological data
  analysis.
\newblock In \emph{Signal and Information Processing (GlobalSIP), 2016 IEEE
  Global Conference on}, pages 1260--1264. IEEE, 2016.

\bibitem[Emrani et~al.(2014)Emrani, Gentimis, and Krim]{emrani2014persistent}
S.~Emrani, T.~Gentimis, and H.~Krim.
\newblock Persistent homology of delay embeddings and its application to wheeze
  detection.
\newblock \emph{IEEE Signal Processing Letters}, 21\penalty0 (4):\penalty0
  459--463, 2014.

\bibitem[Garland et~al.(2016)Garland, Bradley, and Meiss]{garland2016exploring}
J.~Garland, E.~Bradley, and J.~D. Meiss.
\newblock Exploring the topology of dynamical reconstructions.
\newblock \emph{Physica D: Nonlinear Phenomena}, 334:\penalty0 49--59, 2016.

\bibitem[Hoffman et~al.(1998)Hoffman, Kaufman, Halverson, and
  Schrag]{hoffman1998neoproterozoic}
P.~F. Hoffman, A.~J. Kaufman, G.~P. Halverson, and D.~P. Schrag.
\newblock A neoproterozoic snowball earth.
\newblock \emph{science}, 281\penalty0 (5381):\penalty0 1342--1346, 1998.

\bibitem[Khasawneh and Munch(2016)]{khasawneh2016chatter}
F.~A. Khasawneh and E.~Munch.
\newblock Chatter detection in turning using persistent homology.
\newblock \emph{Mechanical Systems and Signal Processing}, 70:\penalty0
  527--541, 2016.

\bibitem[Khasawneh et~al.(2018)Khasawneh, Munch, and
  Perea]{khasawneh2018chatter}
F.~A. Khasawneh, E.~Munch, and J.~A. Perea.
\newblock Chatter classification in turning using machine learning and
  topological data analysis.
\newblock In \emph{14th IFAC Workshop on Time Delay Systems TDS 2018},
  volume~51, pages 195--200. International Federation of Automatic Control,
  2018.

\bibitem[Korn and Faure(2003)]{korn2003there}
H.~Korn and P.~Faure.
\newblock Is there chaos in the brain? ii. experimental evidence and related
  models.
\newblock \emph{Comptes rendus biologies}, 326\penalty0 (9):\penalty0 787--840,
  2003.

\bibitem[Latschev(2001)]{latschev2001vietoris}
J.~Latschev.
\newblock Vietoris-rips complexes of metric spaces near a closed riemannian
  manifold.
\newblock \emph{Archiv der Mathematik}, 77\penalty0 (6):\penalty0 522--528,
  2001.

\bibitem[Lorenz(1963)]{lorenz1963deterministic}
E.~N. Lorenz.
\newblock Deterministic nonperiodic flow.
\newblock \emph{Journal of the atmospheric sciences}, 20\penalty0 (2):\penalty0
  130--141, 1963.

\bibitem[MacPherson and Schweinhart(2012)]{macpherson2012measuring}
R.~MacPherson and B.~Schweinhart.
\newblock Measuring shape with topology.
\newblock \emph{Journal of Mathematical Physics}, 53\penalty0 (7):\penalty0
  073516, 2012.

\bibitem[Perea(2016)]{perea2016persistent}
J.~A. Perea.
\newblock Persistent homology of toroidal sliding window embeddings.
\newblock In \emph{2016 IEEE International Conference on Acoustics, Speech and
  Signal Processing (ICASSP)}, pages 6435--6439. IEEE, 2016.

\bibitem[Perea(2018)]{perea2018brief}
J.~A. Perea.
\newblock A brief history of persistence.
\newblock \emph{preprint arXiv:1809.03624}, 2018.
\newblock \url{https://arxiv.org/abs/1809.03624}.

\bibitem[Perea and Harer(2015)]{perea2015sliding}
J.~A. Perea and J.~Harer.
\newblock Sliding windows and persistence: An application of topological
  methods to signal analysis.
\newblock \emph{Foundations of Computational Mathematics}, 15\penalty0
  (3):\penalty0 799--838, 2015.

\bibitem[Perea et~al.(2015)Perea, Deckard, Haase, and Harer]{perea2015sw1pers}
J.~A. Perea, A.~Deckard, S.~B. Haase, and J.~Harer.
\newblock Sw1pers: Sliding windows and 1-persistence scoring; discovering
  periodicity in gene expression time series data.
\newblock \emph{BMC bioinformatics}, 16\penalty0 (1):\penalty0 257, 2015.

\bibitem[Robins(1999)]{robins1999towards}
V.~Robins.
\newblock Towards computing homology from finite approximations.
\newblock In \emph{Topology proceedings}, volume~24, pages 503--532, 1999.

\bibitem[Ruelle and Takens(1971)]{ruelle1971nature}
D.~Ruelle and F.~Takens.
\newblock On the nature of turbulence.
\newblock \emph{Communications in mathematical physics}, 20\penalty0
  (3):\penalty0 167--192, 1971.

\bibitem[Rustici et~al.(2004)Rustici, Mata, Kivinen, Li{\'o}, Penkett, Burns,
  Hayles, Brazma, Nurse, and B{\"a}hler]{rustici2004periodic}
G.~Rustici, J.~Mata, K.~Kivinen, P.~Li{\'o}, C.~J. Penkett, G.~Burns,
  J.~Hayles, A.~Brazma, P.~Nurse, and J.~B{\"a}hler.
\newblock Periodic gene expression program of the fission yeast cell cycle.
\newblock \emph{Nature genetics}, 36\penalty0 (8):\penalty0 809, 2004.

\bibitem[Takens(1981)]{takens1981detecting}
F.~Takens.
\newblock Detecting strange attractors in turbulence.
\newblock In \emph{Dynamical systems and turbulence, Warwick 1980}, pages
  366--381. Springer, 1981.

\bibitem[Tralie and Berger(2018)]{tralie2018topological}
C.~J. Tralie and M.~Berger.
\newblock Topological eulerian synthesis of slow motion periodic videos.
\newblock In \emph{2018 25th IEEE International Conference on Image Processing
  (ICIP)}, pages 3573--3577, 2018.

\bibitem[Tralie and Perea(2018)]{tralie2018quasi}
C.~J. Tralie and J.~A. Perea.
\newblock (quasi) periodicity quantification in video data, using topology.
\newblock \emph{SIAM Journal on Imaging Sciences}, 11\penalty0 (2):\penalty0
  1049--1077, 2018.

\bibitem[Xu et~al.()Xu, Tralie, Antia, Lin, and Perea]{xu2018twisty}
B.~Xu, C.~J. Tralie, A.~Antia, M.~Lin, and J.~A. Perea.
\newblock Twisty takens: A geometric characterization of good observations on
  dense trajectories.
\newblock \emph{arXiv preprint arXiv:1809.07131}.
\newblock \url{https://arxiv.org/pdf/1809.07131}.

\end{thebibliography}
\end{document}